\documentclass[11pt]{article}
\usepackage{amssymb,amsfonts,amsmath,latexsym,epsf,tikz,url}
\newtheorem{theorem}{Theorem}[section]

\newtheorem{conjecture}[theorem]{Conjecture}
\newtheorem{corollary}[theorem]{Corollary}
\newtheorem{lemma}[theorem]{Lemma}

\newcommand{\proof}{\noindent{\bf Proof.\ }}
\newcommand{\qed}{\hfill $\square$\medskip}

\begin{document}

\title{On the roots of total domination polynomial of graphs}

\author{Saeid Alikhani$^{}$\footnote{Corresponding author}  and  Nasrin Jafari }

\date{\today}

\maketitle

\begin{center}

   Department of Mathematics, Yazd University, 89195-741, Yazd, Iran\\
{\tt alikhani@yazd.ac.ir}\\

\end{center}


\begin{abstract}
Let $G = (V, E)$ be a simple graph of order $n$. The total dominating set of $G$ is a subset $D$ of $V$ that every vertex of $V$ is adjacent to some vertices of $D$. The total domination number of $G$ is equal to minimum cardinality of  total dominating set in $G$ and denoted by $\gamma_t(G)$. The total domination polynomial of $G$ is the polynomial $D_t(G,x)=\sum_{i=\gamma_t(G)}^n d_t(G,i)$, where $d_t(G,i)$ is the number of total dominating sets of $G$ of size $i$. In this paper, we study roots of total domination polynomial of some graphs.  We show that  all roots of $D_t(G, x)$ lie in the circle with center $(-1, 0)$ and the
radius $\sqrt[\delta]{2^n-1}$, where $\delta$ is the minimum degree of $G$. As a consequence we prove that if $\delta\geq \frac{2n}{3}$,
 then every integer root of $D_t(G, x)$ lies in the set $\{-3,-2,-1,0\}$. 

\end{abstract}


\medskip
\noindent{\bf AMS Subj.\ Class.:}  05C69.

\section{Introduction}

Let $G = (V, E)$ be a simple graph. The order of $G$ is the number of vertices of $G$. For any vertex $ v \in V$, the open neighborhood of $v$ is the set $N(v)=\{ u \in V | uv \in E\}$ and the closed neighborhood is the set $N[v]=N (v) \cup \{v\}$.
For a set $S\subset V$, the open neighborhood of $S$ is the set $N(S)=\bigcup_{v\in S }N(v)$ and the closed neighborhood of $S$ is the set $N[S]=N (S) \cup S$. The set $D\subset V$ is a total dominating set if every vertex of $V$ is adjacent to some vertices of $D$, or equivalently, $N(D)=V$. The total domination  number $\gamma_t(G)$ is the minimum cardinality of a total dominating set in $G$. A total dominating set with cardinality $\gamma_t(G)$ is called a $\gamma_t$-set. An $i$-subset of $V$ is a subset of $V$ of cardinality $i$. Let $D_t(G, i)$ be the family of total dominating sets of $G$ which are $i$-subsets and let $d_t(G,i)=|D_t(G, i)|$. The polynomial $D_t(G; x)=\sum_{i=1}^n d_t(G,i)x^i$ is defined as total domination polynomial of $G$. As an example, for all $n\geq 1$ and $2\leq k\leq n$,  $d_t(K_n,i) ={n\choose i}$ and so $D_t(K_n,x) = (x + 1)^n-nx-1$. A root of $D_t(G, x)$ is called a total  domination root of $G$.
For many graph polynomials, their roots have attracted considerable attention. The roots of the
chromatic polynomial, independence polynomial, matching polynomial and reliability
polynomials have been studied extensively \cite{10,11,12,13,14,15}. 
Roots of domination polynomial has investigated in some papers. Graphs with exactly one, two, three and four domination roots characterized in  \cite{AAOP,Euro,four}. Also  Brown and Tufts studied the location of the roots of domination
polynomials for some families of graphs such as bipartite cocktail party graphs and complete bipartite graphs. In particular,
they showed that the set of all domination roots is dense in the complex plane (\cite{Brown}). Oboudi in \cite{Oboudi}, proved that if $\delta$ is
the minimum degree of vertices of $G$ of order $n$, then  all roots of $D(G, x)$ lies in the set
$\{z : |z + 1| \leq \sqrt[\delta+1]{2^n-1}\}$. Motivated by a conjecture
 which states that every integer root of $D(G, x)$ is $-2$ or $0$ (see \cite{Euro}), he has shown that if $\delta \geq \frac{2n}{3}- 1$, then every integer root of $D(G, x)$ is $-2$ or $0$. 
In this paper we would like to study the total domination roots.

The corona of two graphs $G_1$ and $G_2$, as defined by Frucht and Harary in \cite{harary}, is the graph $G= G_1\circ G_2$ formed from one copy of $G_1$ and $|V(G_1)|$ copies of $G_2$, where the i-th vertex of $G_1$ is adjacent to every vertex in the i-th copy of $G_2$. The corona $G\circ K_1$, in particular, is the graph constructed from a copy of $G$, where for each vertex $v \in V(G)$, a new vertex $v'$ and a pendant edge $vv'$ are added. The join of two graphs $G_1$ and $G_2$, denoted by $G_1\vee G_2$, is a graph with vertex set $V(G_1)\cup V(G_2)$ and edge set $E(G_1)\cup E(G_2)\cup \{ uv| u\in V(G_1) ~and~ v \in V(G_2)\}$.

As usual we denote the complete graph, path and cycle of order $n $ by $K_n$, $P_n$ and $C_n$, respectively. Also $K_{1,n}$ is the star graph with $n +1$ vertices.

In the next section we compute the total domination polynomial of some specific graphs. In Section 3, we study the roots of total domination polynomial of a graph. Along the way several problems for further study are also listed.

\section{Total domination polynomial of some graphs }

\begin{figure}
	\begin{center}
		\includegraphics[width=10cm,height=2.3cm]{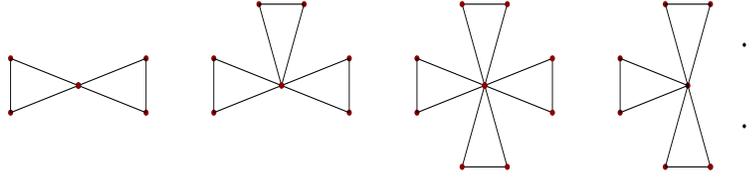}
		\caption{Friendship graphs $F_2, F_3, F_4$ and $F_n$, respectively.}
		\label{Dutch}
	\end{center}
\end{figure}

In this section we shall compute the total domination polynomial of specific graphs. First we consider friendship graph.   
The friendship (or Dutch-Windmill) graph $F_n$ is a graph that can be constructed by coalescence $n$
copies of the cycle graph $C_3$ of length $3$ with a common vertex. The Friendship Theorem of Paul Erd\"{o}s,
Alfred R\'{e}nyi and Vera T. S\'{o}s \cite{erdos}, states that graphs with the property that every two vertices have
exactly one neighbour in common are exactly the friendship graphs.
Figure \ref{Dutch} shows some examples of friendship graphs.

\begin{theorem}\label{tdpf}
	\begin{enumerate}
		\item[(i)] 	For any $n\geq 2$, we have  $\gamma_t(F_n)=2$.
	\item[(ii)] For every $n\geq 2$, $D_t(F_n,x)=x(x+1)^{2n}+x^{2n}-x$.
\end{enumerate}
  \end{theorem}
\proof
\begin{enumerate}
	\item [(i)]
Let $\{v_0,v_1,\ldots ,v_{2n}\}$ be vertex set of $F_n$ and $v_0$ be a common vertex in $F_n$. Then for all $1 \leq i\leq 2n$, $\{v_0, v_i\}$ is a total dominating set for $F_n$, therefore $\gamma_t(F_n)\leq 2$. Since for every graph $G$ we have $\gamma_t(G)\geq 2$, so $\gamma_t(F_n)=2$.
\item[(ii)] 
Suppose that $D$ is a total dominating set of $F_n$ of size $k\leq 2n-1$. Obviously $v_0\in D$ for all total dominating set $D$ of size $k\leq 2n-1$ . To choose $k-1$ other vertices of copies of $C_3$ (except the vertex $v_0$), we have $\binom{2n}{k-1}$ possibilities, so $d_t(F_n,k)={2n\choose k-1}$. For $k\geq 2n$, every subset of vertex set in $F_n$ by size $k$ is a total dominating set for $F_n$. So we have result.\qed
\end{enumerate}

\begin{figure}[!ht]
	\hspace{4cm}
	\includegraphics[width=6.5cm,height=2.3cm]{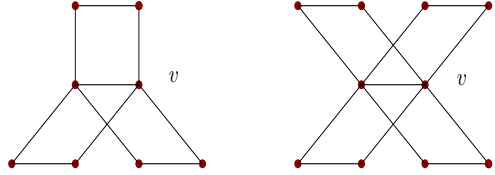}
	\caption{ \label{figure6} The book graphs $B_3$ and $B_4$, respectively.}
\end{figure}

The $n$-book graph $B_n$ can be constructed by bonding $n$ copies of the cycle graph $C_4$ along a common edge $\{u, v\}$, see Figure \ref{figure6}. Here we compute the total domination polynomial of   $n$-book graphs. 

\begin{theorem}\label{B_n}
	\begin{enumerate}
		\item[(i)]
		 	For any $n \geq 2$, $\gamma_t(B_n)=2$.
		 	\item[(ii)]  For every $n\geq 2$,  $D_t(B_n,x)=x^2(x+1)^{2n}+2x^{n+1}(x+1)^n+(2n+1)x^{2n}$. 	
\end{enumerate}
\end{theorem}
\proof
\begin{enumerate}
	\item[(i)] 
	Let $\{u,v\}$ be vertices on the common edge in $B_n$, $\{u,v\}$ is a total dominating set of $B_n$, this implies that $\gamma_t(B_n) \leq 2$. So we have the result.

\item[(ii)] Let $D$ be a total dominating set of $B_n$. Observe that 
 all total dominating set $D$ of size $k\leq n$, contain two vertices $v$ and $u$.
To choose $k-2$ other vertices of copies of $C_4$ (except two vertices $\{v,u\}$), we have $\binom{2n}{k-2}$ possibilities, so for $k\leq n$, $d_t(B_n,k)={2n\choose k-2}$.

For $n+1\leq k\leq 2n-1$, we have two cases: 

Case 1) Suppose that  $\{u,v\}\in D$. To choose $k-2$ other vertices of copies of $C_4\setminus \{v,u\}$'s, we have $\binom{2n}{k-2}$ possibilities.

Case 2) In this case $u \in D$ or $v\in D$. To choose $k-1$ other vertices of $C_4\setminus \{v,u\}$'s, we have $\binom{n}{k-n-1}$ possibilities. So in this case we have $2\binom{n}{k-n-1}$ possibilities.

For $k\geq 2n$, every subset of $V(B_n)$ of  size $k$ is a total dominating set, so for $k\geq 2n$, $d_t(B_n,k)={2n+2\choose k }$. Therefore we have
the  result. \qed
\end{enumerate}

Let $G$ be a graph of order $n$ with vertex set $\{v_1,v_2,\ldots ,v_n\}$ and $H$ be a graph of order $m$ with vertex set $\{u_1,\ldots ,u_m\}$, for every $1\leq i\leq n$,  add $m$ new vertices $\{u^i_1, u^i_2,\ldots,u^i_m \}$ and join $u^i_j$ to $v_i$ for every $1\leq i\leq n$ and $1\leq j\leq m$. By definition this graph is $G\circ H$. The following theorem gives the total domination number of $G\circ H$. 
\begin{theorem}
	For any graph $G$ of order $n$, where $n\geq 2$  and any graph $H$ of order $m$, $\gamma_t(G\circ H)=n$.
\end{theorem}
\proof
	Let $D$ be a total dominating set of $G\circ H$. Since the set  $\{v_1,v_2,\ldots,v_n\}$ is a total dominating set of $G\circ H$, so $|D|\leq n$. Also for every $1\leq i\leq n$,  $v_i\in D$ or at least one of the $u_j^i \in D$, where $1\leq j\leq m$, therefore $|D|\geq n$ and so we have the result.\qed


\begin{theorem}
	For any graph $G$ of order $n\geq 2$ and every $k$, $n\leq k\leq (m+1)n$, we have $d_t(G\circ \overline{K_m},k)=\binom{mn}{k-n}$. Hence $D_t(G\circ \overline{K_m}, x)=x^n(x+1)^{mn}$.
\end{theorem}
\proof
	Suppose that $V(G)=\{v_1,...,v_n\}$ and $D$ is a total dominating set of $G\circ \overline{K_m}$ of size $k\geq n$. Obviously $v_i\in D$ for all $1\leq i\leq n$. To choose $k-n$ other vertices of copies of $K_m$'s, we have $\binom{mn}{k-n}$ possibilities. So we have results.\qed

\begin{theorem}
For every graph  $H$ of order $n$, $D_t(K_1\circ H,x)=x(1+x)^n-x+D_t(H,x)$.
\end{theorem}
\proof
We have two cases
for a total dominating set $D$ of size $k$ in the graph of the form $K_1\circ H$.

Case 1) The set $D$ includes $u$ (the vertex originally in $K_1$) and for choose other vertices in $D$ we have $\binom{n}{k-1}$ possibilities. So the generating function for the number of total dominating sets of graph in this case is $x(1+x)^n-x$.

Case 2) The set $D$ does not include $u$ and it is exactly a total dominating set of $H$. In this case $D_t(H,x)$ is the generating function.

By addition principle, we have $D_t(K_1\circ H,x)=x(1+x)^n-x+D_t(H,x)$.\qed

Since  $D_t(\overline{K_m}\circ H,x)=\prod_{i=1}^m D_t(K_1\circ H,x)=(D_t(K_1\circ H,x)^m$, so we have the following corollary. 
\begin{corollary}
For every graph $H$ of order $n$, we have 
\begin{center}
$D_t(\overline{K_m}\circ H,x)=(x(1+x)^n-x+D_t(H,x))^m$.
\end{center}
\end{corollary}

A finite sequence of real numbers $(a_0, a_1, a_2,...,a_n)$ is said to be unimodal
if there is some $k\in \{0, 1, ..., n\}$, called the mode of sequence, such that
$a_0\leq a_1\leq ...\leq a_{k-1}\leq a_k\geq a_{k+1}\geq ...\geq a_n$;
the mode is unique if $a_{k-1} <a_k> a_{k+1}$. A polynomial is called unimodal
if the sequence of its coefficients is unimodal. The following results prove that 
the total domination polynomial of the complete graph and the friendship graphs are unimodal.

\begin{theorem}
\begin{enumerate} 
\item[(i)] For every $n\in \mathbb{N}$, $D_t(K_n,x)$ is unimodal. 
\item[(ii)] For every $n\geq 2$, $D_t(F_n,x)$ is unimodal. 	
	
\end{enumerate} 	
	\end{theorem} 
\proof\begin{enumerate}
	\item [(i)] 
	Since $D_t(K_n,x)=(x + 1)^n-nx-1={n\choose 2}x^2+{n\choose 3}x^3+...+x^n$ is unimodal, the result follows. 
	\item[(ii)] 
	By Theorem \ref{tdpf},  $D_t(F_n,x)=x(x+1)^{2n}+x^{2n}-x$ and so
	$$D_t(F_n,x)={2n\choose 1}x^2+{2n\choose 2}x^3+....+{2n\choose 2n-2}x^{2n-1}+(2n+1)x^{2n}+x^{2n+1},$$
	since ${2n\choose 2n-2}\geq 2n+1$ is true for all natural number $n\geq 2$, therefore we have the result. \qed	  
\end{enumerate}

Motivated by unimodality of total domination polynomial of specific graphs such as $K_n$ and $F_n$, we have a similar conjecture to unimodality of domination polynomial (\cite{ars}) for total domination polynomial. Notice that that conjecture for domination polynomials is still open. Here, similar to domination polynomial of a graph,  we state and prove the following theorem: 

\begin{theorem}
Let $G$ be a graph of order $n$. Then for every $0\leq i <\frac{n}{2}$, we
have $ d_t(G, i)\leq  d_t(G, i+1)$.
\end{theorem} 
\proof
The proof follows the proof in \cite{ars} with some minor changes.
 Consider a bipartite graph with two partite sets $X$ and $Y$. The
vertices of $X$ are total dominating sets of $G$ of cardinality $i$, and the vertices of
$Y$ are all $(i+1)$-subsets of $V(G)$. Join a vertex $A$ of $X$ to a vertex $B$ of
$Y$, if $A\subseteq B$. Clearly, the degree of each vertex in $X$ is $n-i$. Also for any
$B\in Y$, the degree of $B$ is at most $i+1$. We claim that for any $S\subseteq X$,
$|N(S)|\geq |S|$ and so by Hall's Marriage Theorem, the bipartite graph has a
matching which saturate all vertices of $X$. By contradiction suppose that
there exists $S\subseteq  X$ such that $|N(S)| < |S|$. The number of edges incident
with $S$ is $|S|(n-i)$. Thus by pigeon hole principle, there exists a vertex
$B\in Y$ with degree more than $n-i$. This implies that $i + 1\geq n-i+1$.
Hence $i\geq n$, a contradiction. Thus for every $S\subseteq X$, $|N(S)|\geq |S|$ and the
claim is proved. Since for every $A\in  X$, and every $v\in V(G)\setminus A$, $A\cup \{v\}$ is
a total dominating set of cardinality $i+1$, we conclude that $d_t(G, i+1)\geq d_t(G, i)$
and the proof is complete.\qed 

We end this section with the following conjecture: 
\begin{conjecture}
The total domination polynomial of a graph is unimodal. 
\end{conjecture}

\section{Roots of total  domination polynomial}

In this section, we shall study the roots of the total domination polynomial of a graph. First we state and prove the following theorem: 
\begin{theorem}\label{inclusion}
Let $G=(V,E)$ be a graph of order $n$. Then 
$$D_t(G,x)=\sum_{S\subseteq V}(-1)^{|S|}(x+1)^{n-|N(S)|}.$$
\end{theorem}
\proof
We know  that $(x+1)^n$ is the generating function of the sequence of the number of subsets of $V(G)$. On the other hand
for any $S\subseteq V(G)$, $(x + 1)^{n-|N(S)|}$ is the generating function of the sequence of the subsets of $V(G)$ which has no neighbor in
$S$. By  the inclusion–exclusion principle, we have the result.\qed

As a consequence of Theorem \ref{inclusion}, we have the following theorem about the location of roots of the total domination polynomial. The proofs of Theorems \ref{location}  and Corollary \ref{integer} follows the proof in \cite{Oboudi} with some minor changes:

\begin{theorem} \label{location}
Let $G$ be a graph of order $n$ and  $\delta $ be the minimum degree of $G$.  Then all roots of $D_t(G, x)$ lie in the circle with center $(-1, 0)$ and the
radius $\sqrt[\delta]{2^n-1}$. In other words, if $D_t(G, z) = 0$, then $|z + 1| \leq \sqrt[\delta]{2^n-1}$. 
\end{theorem} 
\proof 
By Theorem \ref{inclusion} and putting   $y=z+1$, we have 
$$D_t(G,z)=\sum_{S\subseteq V}(-1)^{|S|} y^{n-|N(S)|}= y^n-f(y),$$
where $f(y)= \sum_{\emptyset\neq S\subseteq V}(-1)^{|S|} y^{n-|N(S)|}$. If $S$ is a non-empty set of $V(G)$, then $|N(S)|\geq \delta$. Suppose that $|y|\geq 1$, then for 
every $\emptyset\neq S\subseteq V$, $|y|^{n-|N(S)|}\leq |y|^{n-\delta} $, and so $|f(y)|\leq (2^n-1)|y|^{n}$. Therefore we have $|y^n-f(y)|>0$ and so $D_t(G,z)>0$. So we have the result. \qed

\begin{corollary}\label{integer}
         Let $G$ be a graph of order $n$. If $\delta=\delta(G)\geq \frac{2n}{3}$, then every integer root of $D_t(G, x)$ lies in the set $\{-3,-2,-1,0\}$.
\end{corollary} 
\proof
 Since $\delta\geq \frac{2n}{3}$, 
$$\sqrt[\delta]{2^n-1}\leq \sqrt[\delta]{2^n}\leq \sqrt{8}.$$
Let $z$ be a root of $D_t(G, x)$. By Theorem \ref{location}, 
$|z+1|\leq \sqrt[\delta]{2^n-1}$. So $|z+1|\leq \sqrt{8}$ and therefore we have the result.\qed

\begin{figure}[h]‎
	‎\begin{minipage}{7.5cm}‎
		‎\includegraphics[width=\textwidth]{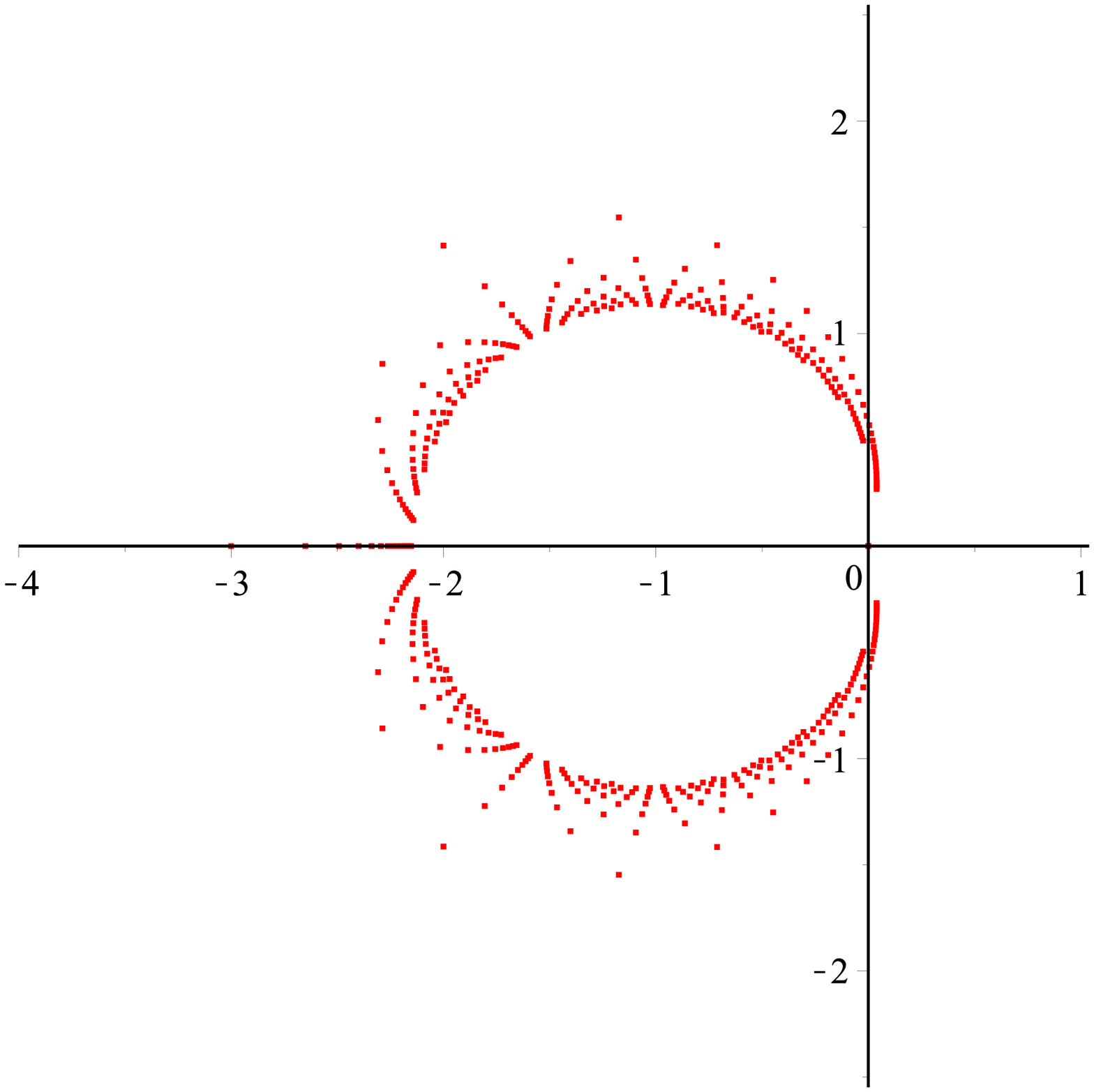}‎
		‎\end{minipage}‎
	‎\begin{minipage}{7.5cm}‎
		‎\includegraphics[width=\textwidth]{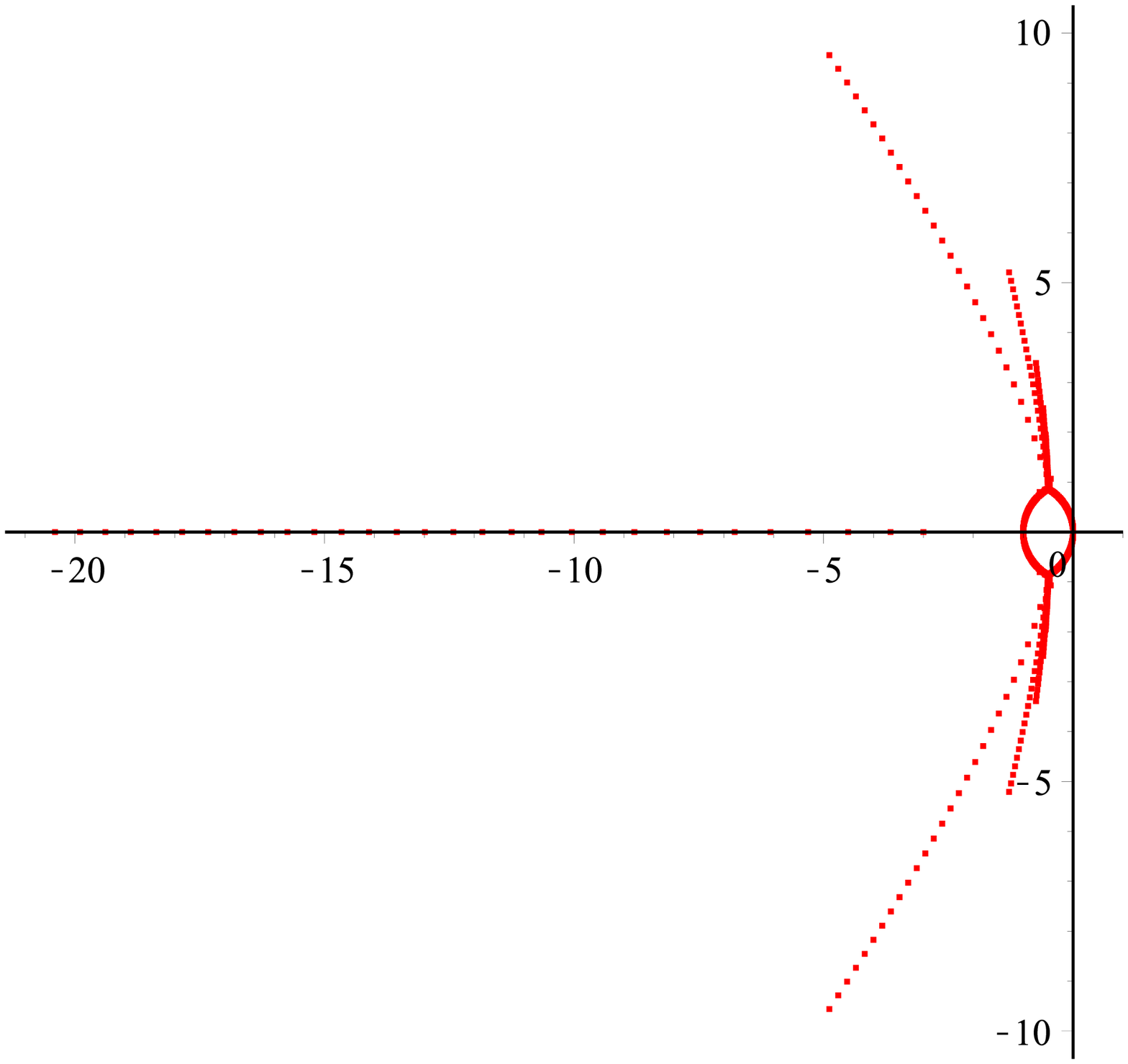}‎
		‎\end{minipage}‎
	‎\caption{\label{figure9} Total domination roots of graphs $K_n$ and $F_n$, for $1 \leq n \leq 20$, respectively.}‎
	‎\end{figure}‎

From Figure \ref{figure9} for complete graph, we state and prove the following theorem: 

\begin{theorem}
       For every even natural number $n$, no nonzero real number is a domination root of $K_n$.
\end{theorem}
\proof
       It is easy to see that  for every $n \in \mathbb{N}~, ~D_t(K_n, x) = (1 + x)^n-nx-1$. If $D_t(K_n, x) = 0$, then for $x\neq 0$, we have 
        \begin{center}
        	$(1 + x)^n = nx + 1$ 
        \end{center}
      We consider two cases, and show in each two there  is no nonzero real solution.

      Case 1) If $x< \frac{-1}{n}$, then  $nx+1<0$, but  for every even number $n\in \mathbb{N}$, $(x+1)^n\geq 0$. It is a contradiction.

      Case 2) If $\frac{-1}{n}\leq x<0$, then we have two following inequalities: 
      \begin{center}
      $(\frac{n-1}{n})^n\leq (x+1)^n<0,~~~~~~~~0<-nx\leq 1$
      \end{center}
      Therefore by adding these inequalities, we have, $(\frac{n-1}{n})^n<(1+x)^n-nx<1$ and therefore cannot be equal to $1$, a contradiction.
              Thus in any event, there are no nonzero real total domination roots of $K_n$ where $n$ is even. \qed

Here we shall investigate the roots of total domination polynomial of the friendship and the book graphs. First we consider the roots of $D_t(F_n,x)$.
Note from Figure \ref{figure9}  it appears that  $D_t(F_n, x)$ has negative real roots of arbitrarily large modulus. This is indeed so, and the following lemma will be used in our proof.

\begin{lemma}{\rm\cite{Brown}} \label{limit}
	$$ lim_{n \rightarrow \infty} ln(n)\Big(\frac{ln(n)-1}{n}\Big)^n=0.$$
\end{lemma}

The basic idea  of the following result follows from  the proof of Theorem $8$ in \cite{Brown}.  

\begin{theorem} 
	The total domination polynomial of the friendship  graph, $D_t(F_n, x)$, 
	where $n\geq2$, 	has a real root in the interval $(-n,-ln(n))$, for $ n$ sufficiently large.
\end{theorem} 
\proof 
Suppose that $$f_{2n}(x)=D_t(F_n,x)=x(1+x)^{2n}+x^{2n}-x.$$
Observe that 
$$f_{2n}(x)=x^{2n+1}+(2n+1)x^{2n}+{2n\choose 2n-2}x^{2n-1}+{2n\choose 2n-3}x^{2n-2}+...+{2n\choose 1}x^2.$$
Consider

$$f_{2n}(-n)=(-1)^{2n+1}n^{2n+1}\Big(1-\frac{2n+1}{n}+\frac{{2n\choose 2}}{(n)^2}-...+\frac{(-1){2n\choose 2n-1}}{(n)^{2n-1}}\Big).$$ 

It is not difficult to see that the following inequality is true  for $n$ sufficiently large, which implies $f_{2n}(-n)<0$  for $n$ sufficiently large,

$$ \frac{2n+1}{n}-\frac{{2n\choose 2}}{(n)^2}+...-\frac{(-1){2n\choose 2n-1}}{(n)^n}<1.$$

Now consider
\begin{eqnarray*}
	f_{2n}(-ln(n))&=&(-ln(n))^{2n}+(-ln(n))\big((1-ln(n))^{2n}-1\big)\\ &=&(ln(n))^{2n}\Big(1-ln(n)\big(\frac{1-ln(n)}{ln(n)}\big)^{2n}+ln(n)\frac{1}{(ln(n))^{2n}}\Big).
\end{eqnarray*}

From Lemma \ref{limit}, we have $ ln(n)\Big(\frac{(ln(n)-1)}{n}\Big)^n \rightarrow 0$, as $n\rightarrow \infty$ which implies that
$f_{2n}(− ln(n))>0$. By the Intermediate Value Theorem, for sufficiently large
$n$, $f_{2n}(x) = D_t(F_n, x)$ has a real root in the interval $(-n,-ln(n))$.\qed

\begin{figure}[!ht]
	\hspace{3cm}
	\includegraphics[width=8cm,height=8cm]{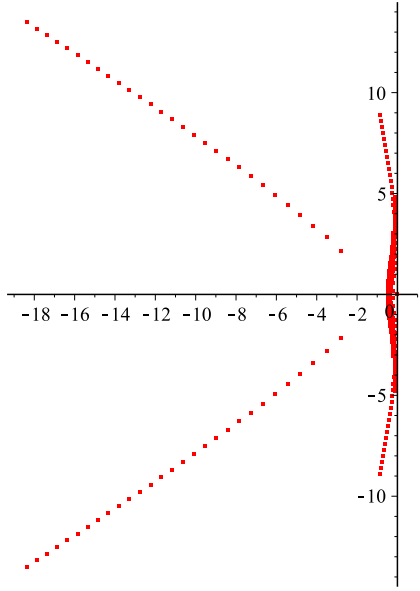}
	\caption{ \label{figure10} Total domination roots of graphs $B_n$ , for $2 \leq n \leq 30$.}
\end{figure}

The following theorem is about the total domination roots of book graph $B_n$ (see Figure \ref{figure10}).

\begin{theorem}
          For every natural number $n$, no nonzero real number is a total domination root of $B_n$
\end{theorem}
 \proof
         By Theorem \ref{B_n}, if $D_t(B_n, x) = 0$, then for  $x\neq 0$ we have
         \begin{center}
         $(x(x+1)^n+x^n)^2=-2nx^{2n}$.
         \end{center}
         Obviously the above equality is true just for real number 0, since for nonzero real number the left side of equality is positive but the right side is negative.\qed

Finally we consider the location of the roots of total domination polynomial of $K_{m,n}$. 

\begin{theorem}
      Let $G$ be the complete bipartite graph$,K_{m,n }= (V_1\cup  V_2,E)$. Then 
     \begin{center}
            $D_t(K_{m,n}, x) = (x+1)^{m+n} - (x+1)^m -(x+1)^n + 1$.
     \end{center}\label{Kmn}
\end{theorem}
\proof
Let $i$ be a natural number $1\leq  i \leq m+n$. We shall determine $d_t(K_{m,n},i)$. Every set $S\subseteq V_1\cup V_2$, with $S\cap V_1\neq \emptyset$ and $S\cap V_2\neq \emptyset$  is a total dominating set of $K_{m,n}$, but subsets which contain only vertices of $V_1$ or $V_2$ are not. So we have the result. \qed

\begin{figure}[!ht]
	\hspace{3cm}
	\includegraphics[width=8cm,height=8cm]{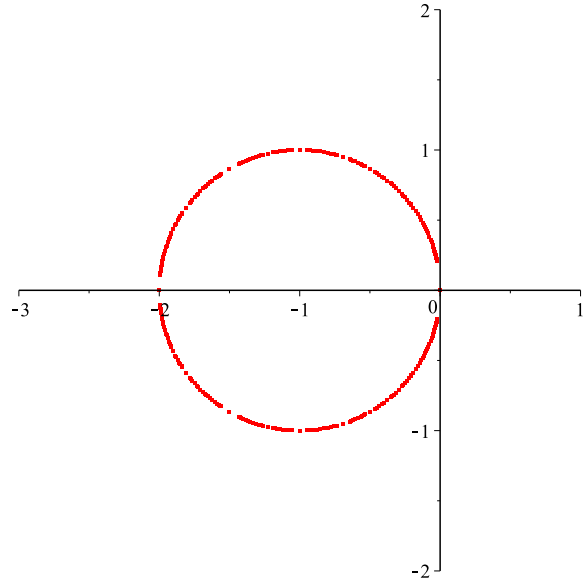}
	\caption{ \label{figure11} Total domination roots of graphs $K_{m,n}$ , for $1 \leq n \leq 30$.}
\end{figure}

Figure \ref{figure11} shows the total domination roots of $K_{m,n}$ for  $1 \leq n \leq 30$. We state and prove the following theorem: 
\begin{theorem} 
	All roots of $D_t(K_{m,n},x)$ are on the circle with center $(-1,0)$ and radius one. 
   \end{theorem}
\proof
    By theorem \ref{Kmn}, for every natural numbers $m,n$, 
    \begin{center}
    $D_t(K_{m,n}, x) = (x+1)^{m+n} - (x+1)^m -(x+1)^n + 1$. 
    \end{center}
    If $D_t(K_{m,n},z)=0$, then we have
     \begin{center}
           $(z+1)^{m+n} - (z+1)^m =(z+1)^n - 1$\\
           and\\
           $(z+1)^{m+n} - (z+1)^n =(z+1)^m - 1$.
    \end{center}
    Therefore $((z+1)^m-1)((z+1)^n-1)=0$, and so we have the result. \qed
    
  We end the paper by proposing the following conjecture: 
    
\begin{conjecture}
Let $G$ be a graph. If $r$ is an integer root of $D_t(G,x)$, then $r\in\{-3,-2,-1,0\}$.
\end{conjecture}


\end{document}